\documentclass[11pt]{amsart}
\usepackage{psfrag,graphicx}
\usepackage{amsmath,amsthm,epsf,amscd,amssymb,verbatim}
\usepackage{palatino}

\newtheorem{thm}{Theorem}[section]
\newtheorem{lem}[thm]{Lemma}
\newtheorem{cor}[thm]{Corollary}
\newtheorem{prop}[thm]{Proposition}

\theoremstyle{definition}

\theoremstyle{remark}

\newtheorem{quest}[thm]{Question}
\numberwithin{equation}{section}
\font\bbb=msbm10 scaled 1100

\newcommand{\be}{\begin{equation}}
\newcommand{\ee}{\end{equation}}
\newcommand{\bea}{\begin{eqnarray}}
\newcommand{\eea}{\end{eqnarray}}
\newcommand{\bmini}{\footnotesize\begin{center}\begin{minipage}{5.5in}}
\newcommand{\emini}{\end{minipage}\end{center}\normalsize}
\newcommand{\pf}{{\em Proof: }}
\newcommand{\R}{\mbox{\bbb R}}

\newcommand{\Z}{\mbox{\bbb Z}} 
\newcommand{\real}{\mbox{\bbb R}}

\newcommand{\nats}{\mbox{\bbb N}}

\newcommand{\zed}{\mbox{\bbb Z}} 
\newcommand{\eg}{{\em e.g.}}

\newcommand{\dbyd}[2]{\frac{\displaystyle\partial #1}{\displaystyle\partial #2}}

\newcommand{\norm}[1]{\left\|{#1}\right\|}
\newcommand{\grad}{\nabla}

\renewcommand{\div}{\nabla\cdot}

\newcommand{\curl}{\nabla\times}

\newcommand{\id}{{\mbox{{\sc Id}}}}

\newcommand{\G}{{\mathcal G}}

\newcommand{\ith}{{i^{th}}}
\newcommand{\M}{{\mathcal M}}
\begin{document}

\title{GENERIC HYDRODYNAMIC INSTABILITY OF CURL EIGENFIELDS}
\author{John Etnyre}
\address{Department of Mathematics, University of Pennsylvania, 
        Philadelphia PA, 19104-6395}
\thanks{JE supported in part by NSF Grant \# DMS - 0203941 .}%%%

\author{Robert Ghrist}
\address{Department of Mathematics, University of Illinois, 
        Urbana IL, 61801}
\thanks{RG supported in part by NSF Grant \# DMS - 0134408.}

\subjclass{Primary: 76E09, 37J55; Secondary: 76B99, 53D40}
\keywords{hydrodynamic instability; Beltrami field; Reeb field; 
  contact homology}

\begin{abstract}
We prove that for generic geometry, the curl-eigenfield
solutions to the steady Euler equations on $\real^3/\zed^3$ 
are all hydrodynamically unstable (linear, $L^2$ norm). 
The proof involves a marriage of contact topological methods 
with the instability criterion of Friedlander-Vishik. An 
application of contact homology is the crucial step.
\end{abstract}

\maketitle

% ************************************************************************
\section{Introduction and statement of results}
% ************************************************************************

The problem of linear hydrodynamic instability
for steady Euler flows on three-dimensional domains
is classical in nature and foundational in implication.
It is universally asserted that in dimension three 
such flows are almost always unstable, though the precise 
definition of {\em ``almost always''} is an issue left
undiscussed. The small literature on generic properties
of fluid flows \cite{FT,ST} focuses on the Navier-Stokes 
setting and uses external forcing or Dirichlet data as
a parameter.

We present a clear formulation of the problem and prove a 
generic instability theorem for a large class of flows --- the 
curl eigenfields --- which form the most fascinating 
and challenging steady solutions to the Euler equations. 
The chief difficulty with genericity issues for curl eigenfields is that
the ``space of all eigenfields'' on a typical Riemannian 
three-manifold is a ``discrete'' space and is not 
amenable to perturbations. Our idea in formulating a 
well-defined genericity statement is to use the geometry of the 
domain as a parameter.

{\bf Main Theorem:} 
{\em 
For a generic set of $C^r$ Riemannian metrics on $T^3:=\real^3/\zed^3$
(for each $2\leq r<\infty$), 
all of the curl-eigenfield solutions to the Euler equations
(with nonzero eigenvalue) are linearly hydrodynamically
unstable in energy norm.  
}

Here and throughout the paper, the term {\em generic} is 
interchangeable with the term {\em residual}: a subset
$A\subset X$ is residual if it is a countable intersection of 
open dense subsets of $X$. All of the function spaces
dealt with in this paper (spaces of vector fields, etc.) are
Baire spaces, implying that residual subsets are dense, 
though not necessarily open. We use the language and notations
of differential forms throughout the paper: $d$ denotes the 
exterior derivative, $\iota$ denotes contraction, $*$ denotes
the Hodge-star operator, $\delta:=*d*$ is a codifferential,
and $\Delta:=d\delta+\delta d$ is the Laplacian.

Given the recent excellent surveys on the intricacies of the 
instability problem \cite{F, FY}, there is little need to 
reintroduce the perspectives in detail. In brief, a vector field 
$u$ is a solution to the {\em steady Euler equations} if
\begin{equation}
\label{eq_euler}
        (u\cdot\nabla)u = -\grad p \,\,\,\, ; \,\,\,\, \div u=0 ,
\end{equation}
for some real-valued pressure function $p$. Given such a solution,
$u$ is said to be {\em linearly stable} if, for every sufficiently 
small divergence-free field $v(0)$, the evolution of $v(t)$ under 
the linearized Euler equation about $u$, 
\begin{equation}
\label{eq_lineuler}
        \dbyd{v}{t} + (u\cdot\nabla)v + (v\cdot\nabla)u = -\grad p ,
\end{equation}
is bounded in some predetermined norm. For the remainder of this 
work, we will, following \cite{F, FV}, use the energy ($L^2$) norm
on vector fields. The solution $u$ is said to be {\em linearly
unstable} if, for some $v(0)$, the solution $v(t)$ has unbounded
growth in the chosen norm. 

Thanks to an insight of Arnold and the analysis of 
Friedlander-Vishik \cite{FV} (who used the technique developed by 
Lifshitz-Hameiri \cite{LH}), 
it is now known that the underlying dynamics of the flowlines 
of the steady solution $u$ can force linear instability. 
In particular, we rely on the following:

{\bf Instability Criterion:} \cite{FV} 
{\em
The presence of a nondegenerate periodic orbit of hyperbolic
(saddle) type in a steady Euler flow induces linear instability
in the energy norm.
}

The outline of the proof of the Main Theorem is as follows:
First, we show that for a generic metric, all of the 
fixed points of all the curl-eigenfields in that metric
are nondegenerate. Any nondegenerate fixed point is immediately
of saddle type (since the flow is divergence-free) and thus forces 
instability. This, then, provides a quick proof of generic
instability for invariant flow domains such as $S^2\times[0,1]$
which are forced to have fixed points on the boundaries. 

Most flow domains, however, have vanishing Euler characteristic
and thus admit flows without fixed points. 
The crux of the difficulty is determining when an eigenfield 
without fixed points possesses a saddle type periodic orbit.
We show that in the absence of fixed points, all of the periodic
orbits are generically nondegenerate. But how does one determine
if periodic orbits exist and are of saddle type? 
To answer this last, most difficult question, we turn to 
contact topological methods, including the recently announced 
{\em contact homology} of Eliashberg, Givental, and Hofer
\cite{EGH}. For an introduction to contact topological techniques
in fluid dynamics, see \cite{EG:I, GK01}. 

We restrict attention in this work to $\real^3$ with periodic
boundary conditions: $T^3$. The sole impediment to applying the
proof to arbitrary three-dimensional domains is the computation 
of the contact homology. For domains with boundary which admit a 
nonvanishing vector field (such as a solid torus), the 
techniques of \cite{EG:II} should suffice to adapt the proofs
to this setting.

The Instability Criterion exists in a slightly more general 
form: a non-periodic orbit having a positive Lyapunov exponent
is sufficiently expanding to push through the analysis. We have
chosen to focus on instances of strict recurrence (fixed points
and periodic orbits). 

Finally, our results would be much improved by changing the genericity 
parameter to be the shape of the boundary (the space of embeddings
of the boundary into a fixed Euclidean space). That this is 
possible is asserted in \cite{Uhlenbeck} and explored in 
\cite{Henry} in the context of eigenfunctions of the Laplace
operator. We do not treat this important case in this paper.

% ************************************************************************
\section{Generic eigenfields}
\label{sec_generic}
% ************************************************************************

This section presents the basic notation and definitions, then 
continues with the principal technical lemmas. 

% ------------------------------------------------------------------------
\subsection{Curl eigenfields}
\label{sec_beltrami}
% ------------------------------------------------------------------------

We begin with the class of curl eigenfield solutions to the Euler
equations on a compact boundaryless three-manifold $M$ with 
Riemannian metric $g$ and volume form $\mu$ 
(assumed to be the volume form derived from $g$ for simplicity --- this 
is not entirely necessary \cite{EG:I}). 
The vector field $u$ is a curl eigenfield of eigenvalue $\lambda$ iff
\begin{equation}
        \lambda \mu(u,\cdot,\cdot) = d(g(u,\cdot)) ,
\end{equation}
or, in forms notation, if $*d\alpha=\lambda\alpha$, where 
$\alpha:=\iota_ug$ is the dual 1-form to $u$ and $*d$ is the 
curl operator on 1-forms. See \cite{AK} for the notation used in 
geometric fluid dynamics on Riemannian manifolds. 

A {\em Beltrami field} on $M$ is defined to be a volume-preserving 
vector field $u$ whose curl is parallel to $u$: in other words, 
$f\iota_u\mu=d(\iota_ug)$ for some map $f:M\to\real$. It is 
a well-known fact that $f$ is an integral for the flow of $u$:
its values are fixed along orbits of $u$. Thus, in the case when 
$f$ is not constant on open sets, the flowlines are
restricted to invariant 2-tori almost everywhere. This is a 
very restrictive assumption, and is in fact not globally possible on 
most (e.g., hyperbolic) 3-manifolds.
For a typical Riemannian manifold (\eg, one without too many symmetries), 
the moduli space of curl eigenfields is discrete (since the 
curl squared is essentially the Laplacian operator and thus 
possesses discrete
spectrum). The class of Beltrami fields can be significantly 
larger in the case where the dynamics is integrable.  

The Instability Criterion requires some expanding dynamics within
the flow, the simplest examples of which are fixed points and periodic
orbits which are {\em nondegenerate} and of saddle-type. A
nondegenerate fixed point is one whose eigenvalues are all 
nonzero. A nondegenerate periodic orbit for a Hamiltonian flow is
defined to be one whose Floquet multipliers (eigenvalues of the 
linearized return map to a cross-section of the orbit) 
are not equal to one.  

% ------------------------------------------------------------------------
\subsection{Reeb fields}
\label{sec_reeb}
% ------------------------------------------------------------------------

Our analysis of curl-eigenfields is strongly rooted in methods derived
from contact geometry. For an introduction to contact geometry, see
\cite[Ch. 8]{Aeb}. In brief, a {\em contact form} on a 3-manifold $M$
is a 1-form $\alpha\in\Omega^1(M)$ such that $\alpha\wedge d\alpha$ 
is nowhere vanishing. A {\em contact structure} on $M$ is a smooth
plane field on $M$ which is (locally) the kernel of a contact 1-form. 
Two fundamental examples of contact structures are the kernel of 
$dz+x\,dy$ on $\real^3$ and also the plane field orthogonal to the 
fibers of the Hopf fibration of the unit $S^3$ in Euclidean $\real^4$.

To every contact form $\alpha$ is associated a unique vector field, 
called the {\em Reeb field}, which captures the geometry of the 
1-form in the directions transverse to the contact structure. The 
Reeb field of $\alpha$, denoted $X$, is defined implicitly 
via the two conditions:
\begin{equation}
\label{eq_Reeb}
d\alpha(X,\cdot) = 0 \quad ; \quad \alpha(X) = 1. 
\end{equation} 
The dynamics of the Reeb field, together with the geometry of the 
contact structure, suffice to reconstruct the contact 1-form. 

Contact geometry enters fluid dynamics via the following results
\cite{EG:I}. For every nonvanishing curl eigenfield $u$ on $(M^3,g)$
with eigenvalue $\lambda\not=0,$ 
the dual 1-form $\alpha:=\iota_ug$ is a contact 1-form since 
$\alpha\wedge d\alpha = \lambda\alpha\wedge*\alpha\neq 0$. Furthermore, 
the Reeb field of $\alpha$ is a rescaling of $u$. Conversely, for each 
contact form $\alpha$, there is a natural adapted Riemannian metric making 
the Reeb field an eigenfield of the curl operator in that metric: 
\begin{equation}
 g(v,w) := \left(\alpha(v)\otimes\alpha(w)\right) + d\alpha(v,Jw) ,
\end{equation}
where $J$ is any almost-complex structure on $\xi=\ker\alpha$ (a 
bundle isomorphism $J:\xi\to\xi$ satisfying $J^2=-\id$) adapted
to $d\alpha$.

The more general version of this correspondence theorem can be 
used to understand the existence and qualitative behavior
of steady solutions to the Euler equations \cite{EG:I,EG:II,EG:III}.
For example, by exploiting the flexibility of contact forms, one 
can construct steady Euler flows on a 3-d Riemannian ball which 
possesses periodic flowlines exhibiting all knot and link types
\cite{EG:III}. 

Our strategy for the remainder of the paper is first to attack 
the case of curl-eigenfields with fixed points using simple genericity
results. The remaining case concerns (rescalings of) Reeb fields for 
contact forms. We will use contact-topological methods in this last case. 
These topological methods also require a certain degree of nondegeneracy,
and thus necessitate genericity statements.

% ------------------------------------------------------------------------
\subsection{Genericity theorems}
\label{sec_lemmas}
% ------------------------------------------------------------------------

The following technical result provides the basis for the 
instability theorem.

\begin{thm}
\label{thm_generic}
For generic choice of $C^r$ metric (for each $2\leq r<\infty$), 
all of the curl-eigenfields on a compact three-manifold $M$ with 
non-zero eigenvalue have all fixed points nondegenerate. In addition, 
if any eigenfield does not possess any fixed points, then all its 
periodic orbits are nondegenerate. 
\end{thm}

The proof will be detailed through a series of lemmas. The first two 
lemmas establish that generically all non-zero eigenvalues are simple 
and eigenfields are transverse to the zero section. We adapt Uhlenbeck's  
techniques \cite{Uhlenbeck} to our situation with two major modifications:
(1) we work with vector fields rather than functions; (2) we use the
curl operator rather than the Laplacian. Moreover, to simplify 
the analysis we consider the dual situation by regarding the curl 
operator as $\curl := *d$ on the space of 1-forms instead of on vector
fields. Denote the space of Riemannian metrics on $M$ by $\G$ and let
\begin{equation}
E_0 := \{(g,\alpha)\in\G\times\Omega^1(M) : \delta\alpha=0\},
\end{equation}
and
\begin{equation}
E := \ker(*d|_{E_0})^\perp.
\end{equation}
Note $E$ is a bundle over $\G$ and the operator $*d$ is a fibrewise 
map. From the Hodge theorem we know $E=\G\times\delta d \Omega^1$ and 
that $*d:E\to E$ is a bundle isomorphism. Now let 
$S:=\{(g,\alpha)\in E : \|\alpha\|_2=1\}$ and consider
\begin{equation}
  \phi:S\times \R \to E
  \quad;\quad 
  \phi(g,\alpha,\lambda) := \left(g, *d\alpha -\lambda\alpha\right) ,
\end{equation}
so that the inverse image of the zero section gives the curl eigenforms. 
This is (fibrewise) an index zero Fredholm operator to which the transversality
theory detailed in \cite{Uhlenbeck} applies. (Though the details are 
not important, as our operator is elliptic, we will use an appropriate 
Sobolev completion on the fibers of $E$ and H\"older norms on $\G$.) 

\begin{lem}
\label{lem_simple}
For each $r\geq 1$, there exists a residual 
set in the space of $C^r$ metrics on a closed $M^3$ such that  
the eigenspaces of the curl operator (with non-zero eigenvalue)
are 1-dimensional and vary smoothly with the metric.
\end{lem}
\pf
We will show that the zero-section $0$ of $E$ is a regular value of $\phi$.
In this case, following \cite{Uhlenbeck}, $Q:=\phi^{-1}(0)$ is
a manifold that fibers over $\G$ with projection $\pi$. 
A $G_\delta$-dense set of metrics will be regular values of 
$\pi$ and, for these values, $Q_g=\pi^{-1}(g)=\phi_g^{-1}(0)$ 
is a 0-dimensional manifold (here $\phi_g=\phi|_{\pi^{-1}(g)}$). 
For each point $(\alpha,\lambda)$ in $Q_g$ we have $*d\alpha=\lambda\alpha$. 
This $\lambda$ is a simple eigenvalue of curl since $0$ is 
a regular value of $\phi$ (cf. \cite[Lemma 2.3]{Uhlenbeck}).
The eigendecompositions vary smoothly since $Q$ is a manifold.
The smoothness condition 
is required for the application of the Sard-Smale theorem: see 
\cite[p. 1061]{Uhlenbeck}.

In order to check that $0$ is a regular value of $\phi$,
we fix $(g,\alpha,\lambda)\in Q$ and consider $\gamma$ orthogonal 
to the image of $D\phi.$  
Thus for any $(h,\beta,s)\in T_{(g,\alpha,\lambda)}(S\times \R)$ we have
\begin{equation}
\left\langle *d \beta-\lambda \beta+s\alpha+
\left[D_{\G}\phi_{(g,\alpha,\lambda)}\right](h), 
 \gamma\right\rangle =0,
\end{equation}
where differentiation along the direction of $\G$ is denoted 
$[D_{\G}\phi_{(g,\alpha,\lambda)}]$ 
(this is denoted $D_2$ in \cite{Uhlenbeck}).
Thus $\gamma$ is orthogonal to $\alpha$ (and so can be thought of as an 
element of $T_\alpha S$). Moreover, 
$\langle \beta, *d\gamma-\lambda \gamma\rangle=0$
for all $\beta\in T_\alpha S.$ Thus taking $\beta=*d\gamma-\lambda \gamma$
we see that $*d\gamma-\lambda \gamma=0.$ (Note 
$\langle *d\gamma-\lambda \gamma, \alpha\rangle
=\langle *d\gamma,\alpha\rangle=\langle \gamma, 
*d\alpha\rangle=\lambda\langle \gamma,\alpha\rangle=0$
so $*d\gamma-\lambda \gamma\in T_\alpha S.$) 
Hence $\gamma$ is an eigenform for $*d.$
Using $\langle [D_{\G}\phi_{(g,\alpha,\lambda)}](h), \gamma\rangle =0$ 
one may easily show that $\gamma$ must vanish in the neighborhood of some
point, then unique continuation implies $\gamma=0.$ (Alternately one 
could show $\gamma$ vanishes away from the zeros of $\alpha$ and then 
use continuity of $\gamma$.)
\qed

For such a metric, then, one can unambiguously designate the 
$\ith$ eigenfield of curl, for $i\in\nats$. 

\begin{cor}
\label{cor_tangents}
Given any continuous $\beta\in E_\alpha$ there is an $h$ such that
$(h,\beta,0)\in  
T_{(g,\alpha, \lambda)}Q.$
\end{cor}
\pf
Choose any tangent vector $(h,\beta,0)\in T_{(g,\alpha,\lambda)}Q$. 
Being constrained to the tangent space implies that
\begin{equation}
\label{eq_tangent}
 *d\beta-\lambda\beta+ \left[D_{\G}\phi_{(g,\alpha,\lambda)}\right](h)=0.
\end{equation}
From Lemma~\ref{lem_simple} we know that $D_{\G}\phi$ is onto continuous 
1-forms (in $E$). Therefore, given any $\beta$, we can choose a tangent 
perturbation $h$ to the space of metrics such that 
$[D_{\G}\phi](h_i)=\lambda\beta_i-*d\beta_i$, 
thus solving Equation~(\ref{eq_tangent}).
\qed

\begin{lem}
\label{lem_fixedpts}
There is a $G_\delta$ dense subset of $C^r$ metrics in $\G$ ($r\geq 2$)
for which all curl eigenfields with non-zero eigenvalues 
have all fixed points nondegenerate.
\end{lem}
\pf
Following \cite{Uhlenbeck}, consider
\begin{equation}
     \psi : Q\times M \to T^*M \quad ; \quad
     \psi(g,\alpha,\lambda, x):=\alpha(x).
\end{equation}
As in the proof of Lemma~\ref{lem_simple}, it suffices to show that 
the zero-section $0$ is a regular value of $\psi$. 
The smoothness condition %($r>\dim(M)+\dim(0\subset T^*M)-\dim(T^*M)=0$) 
is required for the application of the second transversailty theorem
of \cite[p. 1061]{Uhlenbeck} (cf. proof of Theorem 2, p. 1067).
Consider an arbitrary point $(g,\alpha,\lambda)\in Q$ and denote
by $D_Q\psi$ the derivative of $\psi$ at this point along $Q$.  
For $(g,\alpha,\lambda,x)\in Q\times M$, it is clearly seen that 
$[D_Q\psi](h,\beta,0,0)=\beta(x)$. 
Corollary~\ref{cor_tangents} then implies that the 
derivative $D_Q\psi$ is onto.
\qed

\begin{lem}
\label{lem_perorbs}
For each $i\in\zed-\{0\}$ and each positive integer $T$, there exists 
an open dense set of metrics in $\G$ so that, if the $\ith$ eigenfield of
curl has no fixed points, then all of the periodic orbits of 
period less than $T$ are nondegenerate. 
\end{lem}
\pf
The set of generic (in the sense of Lemma~\ref{lem_simple})
metrics in $\G$ which have no fixed points for the $\ith$
eigenfield is open, so work in this set, denoted $\tilde\G^i$.  
Choose any metric $g\in\tilde\G^i$ and let $\tilde{U}_g^i$ denote a
sufficiently small neighborhood of $g$ in $\tilde\G^i$. 
Let $u$ denote the $\ith$ curl eigenfield of $g$, and $O_g$ the slice
of $Q$ through $(g,u,\lambda)$ lying above $\tilde{U}_g^i.$ 
All Beltrami fields near $u$ are in $O_g.$
Since there are no fixed 
points of $u$, the $g$-dual 1-form $\alpha:=\iota_ug$ is a contact 
form whose Reeb field is a rescaling of $u$ by $1/\norm{u}^{2}$.
Note that rescalings are irrelevant since nondegeneracy 
is a topological property of the flowlines of a vector field. 
The map $(g,u,\lambda)\mapsto (\alpha=\iota_u g, \norm{u}^2, J)$ is a
continuous invertible map (where $J$ is the almost complex structure on
the contact structure discussed above) from $O_g$ to an open set $U$. 
Projecting $U$ to the first factor we obtain an open set $U_\alpha$ of 
contact 1-forms near $\alpha$.

We show there is a dense open set in $U_\alpha$ containing 1-forms with
nondegenerate Reeb vector fields. Then the inverse image of this set
will be open and dense in $U$, which leads to an open dense subset of
$\tilde{U}_g^i$ as desired. 

Let $\alpha'$ be a contact 1-form in $U_\alpha.$  
Gray's theorem (see, e.g., 
\cite[p. 169-171]{Aeb}) says that the perturbed contact structure
$\alpha'$ can be deformed through a contact isotopy to the contact 
structure for $\alpha$. Thus, $\alpha'$ can be deformed to a 1-form 
which is a near-identity rescaling of $\alpha$. From the proof of 
Gray's Theorem (using the Moser method in particular), this isotopy is
smooth with respect to $\alpha'$ --- the entire neighborhood of
1-forms near $\alpha$ can be contact-isotoped to near-identity rescalings
of $\alpha$.

Fixing the family of contact isotopies, there is an induced 
continuous map $\Phi$ from $U_\alpha$ to $C^k(M)$ 
which returns this well-defined scaling function for the 
$\ith$ eigenfield of the domain. 
From results of C. Robinson \cite[Thm. 1.B.iv]{Rob} and 
Hofer et al. \cite[Prop. 6.1]{HWZ}, it follows
that there is an open dense set of near-identity 
rescaling functions $f$ such that the Reeb field for $f\alpha$ has 
all periodic orbits of period $<T$ nondegenerate. 

It remains to show that the $\Phi$-inverse image of this open dense
set is an open dense subset of $U_\alpha$. Since $\Phi$ is
continuous, the inverse image is open. One may easily show denseness
by observing that if $\Phi(\alpha')\alpha=\phi^*\alpha'$ does not have the
desired property then there are functions $f_j$ arbitrarily close to
$f=\Phi(\alpha')$ such that $f_j\alpha$ does  have the desired property and 
\begin{equation}
(\phi^{-1})^*(f_j \alpha)= \frac{f_j\circ\phi^{-1}}{f\circ\phi^{-1}}\alpha'
\end{equation}
approximates $\alpha'.$ 
\qed

{\em Proof of Theorem~\ref{thm_generic}:} 
Take the intersection of the residual sets of metrics from 
Lemma~\ref{lem_simple}, Lemma~\ref{lem_fixedpts}, and Lemma~\ref{lem_perorbs}
over all $i$ and $T\in\nats$. This intersection has the desired properties. 
\qed

% ************************************************************************
\section{Contact homology}
\label{sec_homology}
% ************************************************************************

One of the central problems in the topology of contact structures
is the classification problem: given contact structures $\xi$ and 
$\xi'$ on $M$, is there a isotopy of $M$ which takes $\xi$
to $\xi'$? This problem was greatly clarified by dividing the 
set of contact structures into two mutually exclusive types: 
{\em tight} and {\em overtwisted} structures (for definitions, 
see \cite[p. 192]{Aeb}). A theorem of Eliashberg \cite{Eli} states
that the overtwisted structures are classified by the homotopy
type of the plane field, and thus are easily distinguished. This 
is decidedly not the case for the tight structures, whose classification
is a subtle and challenging problem. 

To this end, Eliashberg, Givental, and Hofer have announced a powerful 
new homology theory for contact structures which uses periodic 
orbits of an associated nondegenerate Reeb field as the
chains, a (shifted) Conley-Zehnder index as the grading, and 
pseudoholomorphic curves in the symplectization of the contact
manifold as the mechanism for a boundary operator \cite{EGH}. 
Very recently, the foundations of this {\em contact homology} have
been given rigorous proofs in the thesis of Bourgeois \cite{Bour}.

% ------------------------------------------------------------------------
\subsection{Introduction to contact homology}
% ------------------------------------------------------------------------

Contact homology is an invariant that counts periodic orbits in
a Reeb field for a contact structure $\xi$. Specifically, if
one fixes a contact structure $\xi$ on $M^3$ then one defines the contact
homology as follows. Choose a contact 1-form $\alpha$ for $\xi$ and 
let $X$ be the corresponding Reeb vector field. For a generic
choice of $\alpha$, all periodic orbits of $X$ are non-degenerate \cite{HWZ}.
Let $\mathcal{C}$ be the set of periodic orbits for the flow of $X$.  
To each element $c\in\mathcal{C}$, a grading, $|c|,$ can be assigned  
using a shifted Conley-Zehnder index --- an integer which 
is approximately equal to the number of half-twists 
the linearized flow performs along one cycle. We do not give
a precise definition as the only feature of the grading 
of concern here is:
\begin{lem}
Any nondegenerate orbit $c$ with $|c|$ odd is hyperbolic.
\end{lem} 
For a rigorous definition and an explanation of this fact see, e.g.,  
\cite{Hof93}.

One now defines the graded algebra $\mathcal{A}$ 
as the free super-commutative unital algebra over
$\Z_2$ with generating set $\mathcal{C}$. This algebra
$\mathcal{A}$ will be the chain ``groups'' for contact homology. 
Before defining the boundary map on $\mathcal{A}$ we note one can refine 
the chains as follows.
If $A$ is an element in the first homology of $M$ then one lets 
$\mathcal{C}_A=\{c\in\mathcal{C}| [c]=A \in H_1(M)\}$ and defines 
$\mathcal{A}_A$ to be the sub-algebra of $\mathcal{A}$ generated 
by $\mathcal{C}_A.$

The boundary operator is defined in terms of holomorphic curves in the
symplectization of $(M,\xi).$ The symplectization of $(M,\xi)$ is 
the four-manifold $W=M\times\real$ with symplectic form 
$\omega:=d(e^t\alpha)$, where $t$ denotes the $\real$ coordinate. 
One now equips $W$ with an almost complex structure $J:TW\to TW$ by 
observing that $T_{(x,t)}W=\xi_x\oplus \real_{\left<X\right>}
\oplus\real_{\left<\frac{\partial}{\partial t}\right>}$,
and by defining $J$ on $\xi$ to be any complex structure on $\xi$ 
compatible with $d\alpha\vert_\xi$      
and to send $X$ to $\frac{\partial}{\partial t}$. This choice of $J$ entwines
the dynamics of $X$ with the $t$-direction. 

One next studies holomorphic curves in $W$; that is, maps 
$\varphi:\Sigma\to W$ from a Riemann surface $(\Sigma, j)$ to $W$ such that 
$d\varphi\circ j=J\circ d\varphi$. It is a fact that there are no 
compact Riemann surfaces in $W$; one must introduce punctures 
\cite{Hof93}. 
If $\Sigma$ is a punctured Riemann surface, the energy of $\Sigma$ is 
defined to be $\int_\Sigma\varphi^*(d\alpha)$. If 
$\varphi=(w,h):\Sigma\to M\times\real$ has 
finite energy then some punctures may be removable. Those that are
not removable are guaranteed to possess a neighborhood 
parametrized by $\{(\theta, \tau): \theta\in S^1
\text{ and } \tau\in[0,\infty)\}$ such that $\lim_{\tau\to\infty} h$ 
approaches $\pm\infty$ and $\lim_{\tau\to\infty} w(\theta,\tau)$ 
approaches a parametrization of a periodic orbit $\gamma$ for
$X$. The intuition behind this is that if a surface has finite 
energy, then in the limit as $t\to\pm\infty$, the surface must be
orthogonal to the contact planes, and thus tangent to the 
$(X,\frac{\partial}{\partial t})$ planes. One calls a puncture 
positive (negative) if $h$ approaches $\infty$ ($-\infty$) as 
$\tau\to\infty.$    
For more information on finite energy holomorphic curves and 
their asymptotics see \cite{Audin}.

Given periodic orbits $a, b_1,\ldots, b_k\in\mathcal{C}$, let 
$\M^a_{b_1\ldots b_k}$ denote the set of finite energy 
holomorphic curves in $W$ with one positive puncture asymptotic 
to $a$ and negative punctures asymptotic to 
$b_1, \ldots, b_k,$ modulo holomorphic reparametrization.
Note: since $J$ is $\real$-invariant, there is an $\real$-action on $\M$.
One now defines 
\begin{equation}
\partial a=\sum \left(\#\M^a_{b_1\ldots b_k}/\real\right)b_1\ldots b_k ,
\end{equation}
where the sum is taken over all $b_1,\ldots b_k$ 
such that the dimension of $\M^a_{b_1\ldots b_k}$ is $1$.

\begin{prop}\cite{Bour,EGH}
\label{CH}
The differential $\partial$ lowers grading by 1.
For a generic contact 1-form (and almost complex structure) 
$\partial^2=0$ and the homology of $(\mathcal{A},\partial)$ is
independent of the contact form chosen for $\xi$ (and the almost 
complex structure). 
\end{prop}

The homology of $(\mathcal{A},\partial)$ is called the 
{\em contact homology} of $(M,\xi)$ and
is denoted $CH(M,\xi).$ It is also useful to consider the 
{\em contractible contact homology}
$CH_0(M,\xi)$ whose chain groups are $\mathcal{A}_{[0]}$, 
the algebra generated by contractible periodic orbits. 
The contractible contact homology is also well defined for 
generic contact forms.

It is in general difficult to find all holomorphic curves in $W$ 
and hence to compute the contact homology. It is frequently 
easier to compute the {\em cylindrical contact homology}.
Here one uses $\mathcal{A}_A$ as the chain groups and when one 
defines the boundary map 
one uses only holomorphic cylinders $\Sigma=S^1\times \real$. 

\begin{prop}\cite{Bour,Ust}
\label{CCH}
The cylindrical contact homology $CCH_A(M,\alpha)$ is well-defined 
and independent of the contact form used so long as there are 
no contractible periodic orbits with grading $-1, 0$ or $1$.
\end{prop}

%-------------------------------------------------------------------
\subsection{The 3-torus}

Proving the hydrodynamic instability theorem requires knowing 
the existence of a periodic orbit of saddle-type for all 
non-degenerate Reeb fields on $T^3$. A deep theorem of Hofer
\cite{Hof93} guarantees that any Reeb field for an overtwisted
contact structure possesses a closed orbit of grading $+1$; thus, 
we need merely cover the case of the tight contact structures. 
For $T^3$, these are fortunately classified \cite{Giroux,Kanda}: there
is an infinite family of isomorphism classes represented by
\begin{equation}
\label{eq_TightT3}
  \xi_k := \ker\left(\sin(kz)dx + \cos(kz)dy\right) ,
\end{equation}
for $k\in\zed-\{0\}$.
The following contact homology argument is the crucial step in the 
instability proof:

\begin{lem}
\label{lem_CHofT3}
For a nondegenerate Reeb field associated to any tight contact structure 
on $T^3$, there is always hyperbolic periodic orbit.
\end{lem}

\begin{proof}
We begin with an explicit cylindrical homology computation which 
controls the grading of orbits. For any generic rescaling of the forms 
in Eqn.~(\ref{eq_TightT3}) for which cylindrical contact homology is 
well-defined, Bourgeois \cite{Bour} 
has shown that $CCH_A(T^3, \xi_k)$ is nontrivial in grading $-1$, where 
$A$ is, say, the homology class given by $\{(x,y,z): y=0,z=0\}$. 
Thus there must be hyperbolic periodic orbits whenever the cylindrical 
contact homology is well defined. If the cylindrical contact
homology is not well-defined then there must be contractible 
periodic orbits with grading either $-1$, $0$ or $1$ via 
Proposition~\ref{CCH}. If the gradings are $1$ or $-1$ then there 
exists a contractible hyperbolic periodic orbit. 

In the case that the grading is $0$, we turn to the full contact homology. 
by using the Bott-Morse perturbation technique of \cite{Bour}, one
shows that there are no contractible periodic orbits for nondegenerate
contact forms close to those of Eqn.~(\ref{eq_TightT3}). Hence, the 
contractible contact homology of $(T^3,\xi_k)$ vanishes. 
Since, in the case considered, we must have a contractible periodic orbit
with grading zero, the chain complex for $CH_0(T^3,\xi_k)$ 
possesses an element $a$ with $|a|=0$. Therefore, there must exist
a nontrivial chain with odd grading which prevents a nonzero cycle 
in the contact homology. This implied chain is the desired hyperbolic orbit.
\end{proof}

% ------------------------------------------------------------------------
\subsection{Proof of main theorem}
% ------------------------------------------------------------------------

\begin{thm}
For generic choice of $C^r$ metric ($2\leq r<\infty$), all of
the curl-eigenfields on a three-torus $T^3$ (with nonzero eigenvalue)
are hydrodynamically unstable.
\end{thm}
\pf
First, use Theorem~\ref{thm_generic} to reduce
everything to either nondegenerate fixed points or periodic orbits. 
Given such a field $u$, if 
it possesses a fixed point, then it is immediately of saddle type
due to volume conservation and satisfies the Instability Criterion.
If the field is free of fixed points, then it is (after a suitable
rescaling which preserves the topology of the flowlines) a 
Reeb field for the contact form $\alpha:=\iota_ug$. If the 
contact structure $\xi=\ker\alpha$ is overtwisted, Hofer's theorem
\cite{Hof93} implies the existence of a periodic orbit with 
grading $+1$. The nondegeneracy implies that
the orbit is of saddle type and hence forces hydrodynamic instability.
In the final case where $\xi$ is tight, the contact homology computation
of Lemma~\ref{lem_CHofT3} implies instability.
\qed

% ************************************************************************
\section{Conclusions and questions}
\label{sec_conc}
% ************************************************************************

With the exception of the contact homology computation for the 
class of tight contact structures on $T^3$, the methods used to 
prove generic instability are applicable to any closed
three-manifold, as well as to compact three-manifolds with
invariant boundaries (\eg, the solid torus). On a few sufficiently
simple three-manifolds (spheres, lens spaces), 
it is possible to have nondegenerate curl eigenfields which 
do not possess hyperbolic periodic orbits.
It remains unclear whether these fields are hydrodynamically
unstable.

Besides the curl eigenfields, the only other steady solutions to 
the Euler equations in 3-d are integrable: the flow domain is 
filled almost everywhere with invariant 2-tori \cite{Arn}. 
The only viable candidates for integrable fields which are not 
hydrodynamically unstable are those which possess a great deal 
of symmetry, \eg, the Hopf flow (unit flow along the fibers of 
the Hopf fibration on the unit 3-sphere in Euclidean $\real^4$).
\begin{quest}
Is the Hopf flow on the round $S^3$ linearly unstable?
\end{quest}

% ************************************************************************
% BIBLIOGRAPHY

%\bibliographystyle{plain}

\end{document}